\documentclass[11pt]{amsart}
\usepackage{amssymb, amsmath, amsthm}
\usepackage[latin1]{inputenc}
\usepackage{graphicx}
\usepackage[final]{hyperref}
\usepackage{color}

\usepackage[a4paper, centering]{geometry}
\usepackage{color}
\usepackage{graphicx}

\newtheorem{theorem}{Theorem}
\newtheorem{proposition}[theorem]{Proposition}
\newtheorem{lemma}[theorem]{Lemma}
\newtheorem{corollary}[theorem]{Corollary}
\theoremstyle{definition}

\theoremstyle{remark}
\newtheorem{remark}[theorem]{Remark}

\def\R{\mathbb{R}}
\def\C{\mathbb{C}}
\def\N{\mathbb{N}}
\def\K{\mathcal{K}}

\def\C{\mathcal{C}}
\def\vps{\varepsilon}

\definecolor{verde}{RGB}{20,150,100}

\newcommand{\EEE}{\color{black}}

\def\conv{\it Conv}

\def\l{\lambda_1}

\def \d{\delta}
\def \e{\varepsilon}

\begin{document}

\title[]{On the honeycomb conjecture \\ for Robin Laplacian eigenvalues }

\author[]{Dorin Bucur, Ilaria Fragal\`a}
\thanks{}

\address[Dorin Bucur]{
Laboratoire de Math\'ematiques UMR 5127 \\
Universit\'e de Savoie,  Campus Scientifique \\
73376 Le-Bourget-Du-Lac (France)
}
\email{dorin.bucur@univ-savoie.fr}

\address[Ilaria Fragal\`a]{
Dipartimento di Matematica \\ Politecnico  di Milano \\
Piazza Leonardo da Vinci, 32 \\
20133 Milano (Italy)
}
\email{ilaria.fragala@polimi.it}

\keywords{ Optimal partitions, honeycomb, Robin Laplacian eigenvalues, $\alpha$-Cheeger constant. }
\subjclass[2010]{52C20, 51M16, 65N25, 49Q10. }

\begin{abstract} We prove that the optimal cluster problem for the sum of the first Robin eigenvalue of the Laplacian,  in the limit of a large number of convex cells,  is asymptotically solved by (the Cheeger sets of) the honeycomb of regular hexagons.  The same result is established for the Robin torsional rigidity.  
  \end{abstract}
\maketitle

\section{Introduction and statement of the results}

Given an open bounded Lipschitz domain  $\Omega$ in $\R ^2$ and a real parameter $\beta \neq 0$, 
 we denote by $\lambda _1 (\Omega, \beta)$  and $\tau (\Omega, \beta)$ the 
first Robin eigenvalue of the Laplacian in $\Omega$ and the Robin torsional rigidity of $\Omega$  with coefficient $\beta$.  They are  defined as
\begin{equation}\label{bfr001}
\l (\Omega, \beta):= \hskip -.5cm \min _{u \in H ^ 1 (\Omega)\setminus \{ 0 \} } \!\!\frac{ \int_{\Omega} |\nabla u | ^ 2 + \beta \int_{\partial \Omega } u ^ 2 }
 { \int _{\Omega} u ^ 2}, \quad \tau ^ {-1} (\Omega, \beta) :=  \hskip -.5cm  \min _{u \in H ^ 1 (\Omega)\setminus \{ 0 \} }\!\! \frac{ \int_{\Omega} |\nabla u | ^ 2 + \beta \int_{\partial \Omega } u ^ 2 }
 { \Big (\int _{\Omega} |u |  \Big ) ^2}.
 \end{equation}
For the eigenvalue problem, the  corresponding Euler-Lagrange equation is  given  by  $$
\begin{cases}
- \Delta u = \l (\Omega, \beta) u & \text{ in } \Omega
\\
\frac{\partial u }{\partial \nu} + \beta u = 0 & \text{ on } \partial \Omega\,.
\end{cases}
$$
For the torsional rigidity, the Euler-Lagrange equation requires more attention  (see for instance \cite{bw17}), specifically in the case $\beta <0$.  For  positive $\beta$, the minimizer solves
$$
 \begin{cases}
- \Delta u = 1 & \text{ in } \Omega
\\
\frac{\partial u }{\partial \nu} + \beta u = 0 & \text{ on } \partial \Omega\,
\end{cases}
$$ 
while for negative $\beta$, the Euler-Lagrange equation may involve a free boundary problem. It is not the purpose of the present paper to discuss this issue, as we focus only on the energy values defined in \eqref{bfr001}.  
Without any attempt of completeness, we refer to  \cite{BG15,BG16,CK16}   for some recent papers in shape optimization involving free boundaries with Robin conditions.

While there is a wide literature about optimal partitions for the first Dirichlet Laplacian eigenvalue (see for instance \cite{BoVe, BNHV, BuBuHe, CTV03, CTV05, He07, He10, HeHoTe, Ramos}), 
to the best of our knowledge the study of the same kind of problem for the first Robin Laplacian eigenvalue
is a completely unexplored field.

Object of this paper are the optimization problems 
\begin{eqnarray}\label{def:p}
&  \displaystyle r _{k} (\Omega, \beta) = \left\{
\begin{array}{ll}
 \inf  \Big \{  \sum_{i = 1} ^ k  \lambda (E _i, \beta) \ :\  \{E_i\} \in \C _k (\Omega) \Big \} \mbox{   if } \beta >0  & \label{f:sum}      \\
 \noalign{\medskip}
\sup  \Big \{  \sum_{i = 1} ^ k  \lambda (E _i, \beta) \ :\  \{E_i\} \in \C _k (\Omega) \Big \} \mbox{   if } \beta <0 
\end{array}
\right.
\, 
\end{eqnarray} 
where
${\C} _k (\Omega)$ denotes the class of convex $k$-clusters $\Omega\subset \R ^2$, meant as families
of $k$ convex bodies contained into $\Omega$ and having mutually disjoint interiors, and, for any $\beta \in \R \setminus \{ 0 \}$, $\lambda (\Omega, \beta)$
may be either  $\lambda _1 (\Omega, \beta)$ or  $\tau^ {-1} (\Omega, \beta) $.

\EEE

We are interested in particular in the asymptotic behaviour of $r _k (\Omega, \beta)$  in the limit as $k \to + \infty$.
Our main motivation is a conjecture due to Caffarelli and Lin  \cite{CaffLin} which predicts that, for the analogous problems in which $\lambda (\Omega, \beta)$ is replaced by the first Dirichlet Laplacian eigenvalue, an optimal configuration is asymptotically given by a packing of regular hexagons,  similarly  to the case of perimeter minimizing partitions settled by Hales  in the celebrated paper \cite{Hales}  (see also \cite{CarMag} for a quantitative formulation). 

Very recently, in \cite{bfvv17} this conjecture has been proved to hold if one takes the Cheeger constant in place of the Dirichlet eigenvalue, and the cells of the partitions are a priori assumed to be convex.  Recall that the  Cheeger constant of $\Omega$ 
(about which a detailed account can be found for instance in \cite{Leo, Pa})
is defined by
\begin{equation}\label{f:defh}
h (\Omega):=\inf \left \{ \frac{{\rm Per} (E, \R ^2)}{|E|}\ :\ E \hbox{ measurable}\, , \ {E\subseteq   \Omega} \right \}\,,
\end{equation}
where ${\rm Per}(  E, \R ^2)$ denotes the perimeter of $E$ in the sense of De Giorgi. 

Clearly, the fact that the notion of Cheeger constant is purely geometrical makes the analysis of optimal partitions, started by Caroccia in \cite{Car17}, much more manageable with respect to the case of eigenvalues. 
Nevertheless, the approach proposed in \cite{bfvv17} does not rely specifically on the definition of $h (\Omega)$, but rather on the validity of a 
discrete Faber-Krahn inequality (which in case of the Cheeger constant has been proved in \cite{BF16}) and on some other geometric properties, like the monotonicity upon inclusions of sets and a scaling behaviour; and in fact, the proof given in \cite{bfvv17} adapts also to other shape functionals for which a polygonal version of Faber-Krahn inequality is available, such as for instance  a power of perimeter or the logarithmic capacity 
(see \cite{SoZa}). 

Now, a polygonal Faber-Krahn inequality for  eigenvalues of the Laplacian  is a long-standing conjecture by P\'olya, for which a proof is still missing (see for instance \cite{H06}); thus
the conjecture by Caffarelli-Lin remains open.

Under Robin boundary conditions, neither for the first Laplacian eigenvalue nor for torsional rigidity, polygonal isoperimetric inequalites are known (even for triangles); furthermore, neither $\lambda _1 (\Omega, \beta)$ nor $\tau (\Omega, \beta)$ behave monotonically under inclusions.   At this point, our results about the honeycomb conjecture for such Robin functionals  
should sound somewhat unexpected. We stress that 
we keep the assumption that the cells of the partitions are convex. In case of the first Robin eigenvalue, we prove: 
\begin{theorem}\label{t:honeycomb}
Let $r _k (\Omega, \beta)$ be defined by \eqref{def:p}, with $\lambda (\Omega, \beta) := \lambda _1 (\Omega, \beta)$. 
Then there holds
$$
\lim _{k\to + \infty} \frac{|\Omega|^ {1/2} }{ k ^ {3/2} } r _k (\Omega, \beta) 
=  \beta h ( H)
\, , 
$$
\smallskip
where $h ( H)$ denotes the Cheeger constant of the unit area regular hexagon. 
\end{theorem}

Theorem \ref{t:honeycomb} is obtained as a consequence of the analogous result proved in \cite{bfvv17} for the Cheeger constant, 
combined with a tight control of the Robin eigenvalues  in terms of the quotient perimeter over area ({\it cf.} Proposition 
\ref{p:bounds} below) when the number of cells is increasing. The idea is that, when there is a great amount of cells $E_i$,
thanks to the non-scale invariance of the Robin eigenvalue, 
for a sufficiently large number of indices $i$, the value of $\lambda _1 (E_i, \beta)$ turns out to be comparable to $\beta |\partial E _i| / |E_i|$, so that the partition behaves like a Cheeger one. From a technical point of view, the key point is to prove that this comparison can be made uniform, except for a negligible number of cells, which do not affect the asymptotical behavior as $k \to + \infty$.

In case of the Robin torsional rigidity, we prove: 

\begin{theorem}\label{t:honeycomb2}
Let $r _k (\Omega, \beta)$ be defined by \eqref{def:p}, with $\lambda (\Omega, \beta) := \tau  ^ {-1}(\Omega, \beta) $. 
Then there holds
$$ \lim _{k\to + \infty} \frac{|\Omega|^ {1/2} }{ k ^ {3/2} } r _k (\Omega, \beta) 
=  \beta h _2( H)
\, , 
$$
\smallskip
where $h_2 ( H)$ denotes the $2$-Cheeger constant of the unit area regular hexagon. 
\end{theorem}

The notion of $2$-Cheeger constant appearing in the statement of Theorem \ref{t:honeycomb2} is a variant of the classical definition \eqref{f:defh} of Cheeger constant; 
precisely, the $2$-Cheeger constant  of a set $\Omega$ is given by
\begin{equation}\label{f:defh2}
h_2 (\Omega):=\inf \left \{ \frac{{\rm Per} (E, \R ^2)}{|E| ^2}\ :\ E \hbox{ measurable}\, , \ {E\subseteq   \Omega} \right \}\,. 
\end{equation}
This generalization of Cheeger constant has already appeared in the literature, actually with the square of volume replaced by  an arbitrary power with exponent $\alpha > 1/2$, see \cite{FiMP, FuMP, PS}. 

In the same fashion as Theorem \ref{t:honeycomb} is obtained by applying the analogous result proved in \cite{bfvv17} for the Cheeger constant, 
combined with a tight control of the Robin eigenvalues  in terms of the quotient perimeter over area, Theorem \ref{t:honeycomb2} is obtained by applying the analogous result for the $2$-Cheeger constant, combined with a tight control of the Robin torsion in terms of the quotient perimeter over the {\it square} of the area. 
Actually, in order to prove Theorem \ref{t:honeycomb}, we need as a first step to  settle a honeycomb-type result for the $2$-Cheeger constant analogous 
to the one proved in \cite{bfvv17} for the Cheeger constant. In turn, this requires to obtain a discrete Faber-Krahn inequality for the $2$-Cheeger constant 
in the vein of \cite{BF16} (but dealing just with convex polygons).

As a consequence of Theorems \ref{t:honeycomb} and \ref{t:honeycomb2}, we can also determine the asymptotic behaviour of similar problems where the energy is of  supremal rather than additive type. Setting
\begin{eqnarray}\label{def:P}
& \displaystyle R _{k} (\Omega, \beta) = \left\{
\begin{array}{ll}
\inf  \Big \{  \max_{i = 1, \dots, k}  \limits \lambda (E _i, \beta) \ :\  \{E_i\} \in \C _k (\Omega) \Big \} \mbox{   if } \beta >0 
& \label{f:max} \\ \noalign{\medskip}
\sup  \Big \{  \min_{i = 1, \dots, k} \limits  \lambda (E _i, \beta) \ :\  \{E_i\} \in \C _k (\Omega) \Big \} \mbox{   if } \beta <0 \,.
\end{array}
\right.
\, 
\end{eqnarray} 
we have: 

\begin{corollary}\label{c:cormax}

(i) If $R _k (\Omega, \beta)$ is defined by \eqref{def:P}, with $\lambda (\Omega, \beta) := \lambda _1 (\Omega, \beta)$, 
there holds
$$
\lim _{k\to + \infty} \frac{|\Omega|^ {1/2}}{ k ^ {1/2}} R _k (\Omega, \beta) 
=\beta h ( H) 
\, ; 
$$

(ii) If $R _k (\Omega, \beta)$ is defined by \eqref{def:P}, with $\lambda (\Omega, \beta) := \tau  ^ {-1}(\Omega, \beta) $, 
there holds
$$ \lim _{k\to + \infty} \frac{|\Omega|^ {1/2}}{ k ^ {1/2}} R _k (\Omega, \beta) 
=
\beta h _2 ( H) 
\, . 
$$

\end{corollary}

The detailed proofs of our results are presented hereafter with the following outline. 

In Section \ref{sec:2C}, we provide some results about the $2$-Cheeger constant which may have an autonomous interest, 
by showing by particular that it satisfies the honeycomb conjecture (with convex cells).

In Section \ref{sec:prel}, we establish some intermediate results towards the proofs of Theorems \ref{t:honeycomb} and \ref{t:honeycomb2}, which are crucial to 
make the connection between optimal Robin partitions and optimal Cheeger partitions: they consist essentially in settling
good upper and lower bounds for $\lambda _ 1 (\Omega, \beta)$ and $\tau(\Omega, \beta)$ in terms of geometrical quantities,  holding under the assumption that the domain $\Omega$ is convex. Actually, in case of
the Robin Laplacian eigenvalue and for $\beta >0$, 
an alternative simpler proof can be obtained   by exploiting for the lower bound more rough  inequalities not requiring convexity; nevertheless, 
since this direct approach is somehow related to the first eigenvalue and does not work for the Robin torsional rigidity, we preferred to follow the same guideline for both cases (see the final Remark \ref{finrem} for more detailed comments in this direction). 

Sections  \ref{sec:proof}, \ref{sec:proof2}, and \ref{sec:proof3} contain respectively the proofs of Theorems \ref{t:honeycomb}, Theorem \ref{t:honeycomb2} and Corollary \ref{c:cormax}.

\section{An auxiliary result about optimal $2$-Cheeger partitions}\label{sec:2C}
The $2$-Cheeger constant shares many features with the classical one. For instance, it is easy to check, by using the same arguments as for the classical Cheeger constant, that a $2$-Cheeger set $C_2(\Omega)$ (namely a solution to problem \eqref{f:defh2}) always exists; moreover a $2$-Cheeger set  is connected, 
its boundary is of class $\mathcal C ^1$ and meets necessarily $\partial \Omega$ (and this occurs tangentially), and 
$\partial C _2 (\Omega ) \cap \Omega$ is made by arcs of circle
(of radius $\big (  2|C_2(\Omega)| h _ 2 (\Omega)| \big ) ^ {-1}$).  
We refer the interested reader to the proofs given in \cite{PS} for a similar notion of $\alpha$-Cheeger constant. 

In this section we present some focused results about the $2$-Cheeger constant, with the final scope of proving that
it satisfies the honeycomb conjecture under convexity constraint on the cells, 
thus extending the result proved in \cite{bfvv17} for the classical Cheeger constant.

We consider the following shape optimization problem, where $\mathcal P _n$ denotes the class of convex polygons with at most $n$ sides:
\begin{equation}\label{p:polygons}\min \Big \{ |\Omega | ^ {3/2} h _2 (\Omega) \ :\ \Omega \in \mathcal P _n \Big \}\,.
\end{equation}
Notice that, since  the functional $\Omega \mapsto h _ 2 (\Omega)$ is  homogenoues  of degree  $-3$ under dilations, multiplying it by 
 $|\Omega | ^ {3/2}$ we find a scale invariant functional. In particular, it turns out that problem \eqref{p:polygons} is well-posed. Moreover, 
in the next lemma we establish that the $2$-Cheeger radius 
 of an optimal polygon is uniquely determined through an explicit representation formula involving just three geometrical quantities: $|\partial \Omega|$, $|\Omega|$, 
 and the functional $\Lambda (\Omega)$ defined by 
\begin{equation}\label{f:Lambda}
\Lambda (\Omega):= \sum _i \cot\Big ( \frac{\theta_i}{2}  \Big ) \, . 
\end{equation}
being  $\theta _1, \dots, \theta _N$ the inner angles of $\Omega$. 
The analogous result for the Cheeger constant 
can be found in \cite{BF16},  in a more general setting not requiring convexity.

 \begin{lemma}\label{l:rep}
There exists a solution to 
problem \eqref{p:polygons}, and any optimal polygon  $\Omega$ admits
 a unique $2$-Cheeger set  $C_2(\Omega)$, which is determined by the equality
\begin{equation}\label{f:charC}
\partial C_2 (\Omega) \cap \Omega = \bigcup \Big \{ \Gamma _\alpha  \ :\  \alpha \in\Theta(\Omega) \Big \}\, ,
\end{equation}
where $\Theta$ denotes the family of inner angles of $\Omega$ and, for any $\alpha \in \Theta$, $\Gamma _\alpha$ is an
arc of circle tangent
 to  the two sides of $\partial \Omega$ forming the angle $\alpha$, of radius
 \begin{equation}\label{f:charr} 
 r _2 (\Omega)=  \frac{|\Omega |} {|\partial \Omega| + \sqrt {|\partial \Omega| ^ 2 - 3 |\Omega| \big (\Lambda (\Omega) - \pi  \big )}}\,.
 \end{equation}
 Moreover, the $2$-Cheeger constant of $\Omega$ is given by
\begin{equation}\label{f:charh}  h _2 (\Omega) =  
\frac{|\partial \Omega| - 2 r _2 (\Omega) \big (\Lambda (\Omega) - \pi  \big )} 
{\big [ |\Omega| -  (r _2 (\Omega)) ^ 2 \big (\Lambda (\Omega) - \pi  \big ) \big ] ^ 2}  
\,.\end{equation}
\end{lemma}

 \proof  The existence of an optimal polygon is straightforward: since we minimize over a closed subclass of the class of  convex polygons a continuous and 
 dilation invariant functionals, it is enough to apply the  direct method of the Calculus of Variations working with the Hausdorff convergence. 
 
In order to get the equality \eqref{f:charC}, assume that 
$\Omega \in {{\mathcal P} _N}$ is a solution to problem \eqref{p:polygons}, and let $C _2(\Omega)$ be a $2$-Cheeger set of $\Omega$.  Then it is readily seen that 
$C_2(\Omega)$ must touch every side of $\Omega$.
(Namely, if by contradiction there exists a side which is not touched by $C_2(\Omega)$, we could construct a domain $\widetilde \Omega$, still belonging to ${{\mathcal P} _N}$, such that
that $C_2(\Omega) \subset \widetilde \Omega \subset \Omega$. Then $|\widetilde \Omega|  < |\Omega|$  and  $ h_2 (\widetilde \Omega) = h_2 (\Omega)$,  so that  
$|\Omega|^ {3/2} h_2 (\Omega) >  |\widetilde \Omega|^ {3/2}    h_2 (\widetilde \Omega)$, against the optimality of $\Omega$). 

As a consequence of the facts that $C_2 (\Omega)$ meets every side of $\Omega$ and it is connected,  we obtain that all the arcs of circle contained into $\partial C_2(\Omega) \cap \Omega$ must be  tangent to two consecutive sides of $\Omega$. 
 
Let us show that, for any $\alpha \in \Theta (\Omega)$,  there exists an arc of type $\Gamma _\alpha$ such that
$\Gamma _\alpha \subseteq \partial C _2(\Omega) \cap \Omega$.  
Let $\alpha \in \Theta (\Omega)$ be fixed, and let $\Omega _{\alpha, r}$ be the domain obtained by ``smoothing'' the corner $\alpha$ by means of an arc of cirumference of radius $r$, tangent to the two sides of $\partial \Omega$ forming the angle $\alpha$. It is readily seen by geometric arguments that, for $r$ sufficiently small, 
$$ {{\rm Per}} ( \Omega _{\alpha, r}, \Omega) = |\partial \Omega| - 2r \cot \big (\frac{\alpha}{2}\big ) + ( \pi - \alpha) r $$
and
$$|\Omega _{\alpha, r}| = |\Omega| - r ^ 2 \cot \big  (\frac{\alpha}{2} \big ) +   \big ( \frac { \pi - \alpha}{2}\big ) r ^ 2  \,.$$
Then, 
$$\frac{{\rm Per} (\Omega _{\alpha, r})}{|\Omega _{\alpha, r}| ^2} =  \frac { |\partial \Omega| - 2r \big [ \tan \big (\frac{\pi- \alpha}{2}\big ) - (\frac{ \pi - \alpha}{2})  \big]   }
{ \Big \{ |\Omega| - r ^ 2 \big [\tan \big  (\frac{\pi- \alpha}{2} \big ) -   \big ( \frac { \pi - \alpha}{2}\big ) \big] \Big \} ^ 2 } \,. $$ 
Since the term in squared parenthesis is positive, we immediately see that the inequality 
$\frac{ {{\rm Per}} ( \Omega _{\alpha, r})  }{|\Omega _{\alpha, r}| ^2} < \frac{|\partial \Omega |  }{|\Omega|^2}$ is satisfied for $r$ sufficiently small. 

\smallskip
In order to get the optimal radius, we have to minimize the function
$$f(r):= \frac { |\partial \Omega| - 2r  \Big [ \sum \limits _{\alpha \in \Theta  (\Omega)} \big [ \tan \big (\frac{\pi- \alpha}{2}\big )   \big]   - \frac{\pi- \alpha}{2} \Big ] }
{\Big \{ |\Omega| - r ^ 2  \Big [\sum \limits_{\alpha \in \Theta (\Omega)} \big [\tan \big  (\frac{\pi- \alpha}{2} \big )  \big]  - \frac{\pi- \alpha}{2}  \Big ] \Big \} ^ 2} =   \frac { |\partial \Omega| - 2
r  \big(  \Lambda (\Omega) - \pi \big )    }{\Big \{|\Omega| - r ^ 2 \big(  \Lambda (\Omega) - \pi \big )   \Big \} ^ 2} \,.$$
 In a neigbourhood of an optimal radius, this function is equal to or larger than the $2$-Cheeger constant of $\Omega$.
In particular, we point out that $f (r)$ is strictly larger than $h _ 2 (\Omega)$  if $r$ is above 
the critical value $\overline r$ for which two distinct arcs, each one tangent to two sides of $\partial \Omega$, lie on the same circumference (this can be easily seen by 
taking as a test in the definition of $h _ 2 (\Omega)$ the  intersection of the two sets obtained by smoothing two consecutive angles   
by an arc of circumference of radius $\overline r + \e$: the quotient between perimeter and squared area of this test is strictly smaller than $f ( \overline r + \e)$). 

Therefore, we proceed to determine the critical points of $f$. 
 
By studying the first derivative $f' (r)$, it is easy to see that $f$ is increasing in the interval $( r_-, r _+)$ between the two roots  of $f'$, 
$$r _\pm :=  \frac{|\Omega|} {|\partial \Omega| \mp \sqrt{ |\partial \Omega| ^ 2 - 3 |\Omega| \big(  \Lambda (\Omega) - \pi \big )  } }\,,$$
so that  $f$ attains its minimum at $r _-$, which gives the value of the $2$-Cheeger radius $r _ 2 (\Omega)$. 
We have thus concluded the proof of \eqref{f:charC}-\eqref{f:charr}, and \eqref{f:charh} follows by definition.

\qed

\bigskip
Relying on  Lemma \ref{l:rep}, we obtain that the $2$-Cheeger constant satisfies the following discrete Faber-Krahn inequality on convex polygons: 

\begin{proposition}\label{p:discreteFK}
The regular $n$-gon solves problem \eqref{p:polygons}, with corresponding infimum equal to
\begin{equation}\label{f:gamma}
\gamma (n):=
\frac{\left(2 \sqrt{n \tan \left(\frac{\pi
   }{n}\right)}+\sqrt{n \tan \left(\frac{\pi }{n}\right)+3
   \pi }\right)^3}{8 \left(n \tan \left(\frac{\pi
   }{n}\right)+\sqrt{n \tan \left(\frac{\pi }{n}\right)}
   \sqrt{n \tan \left(\frac{\pi }{n}\right)+3 \pi }+\pi
   \right)}\,.
   \end{equation}

\end{proposition} 
 \proof
 Thanks to the Lemma \ref{l:rep}, if $\Omega$ is an optimal polygon for problem \eqref{p:polygons}, we can write the cost functional
 $h_2 (\Omega)  |\Omega| ^ {3/2}$ as a function depending only on the isoperimetric quotient $I (\Omega):= \frac{|\partial \Omega| }{|\Omega| ^ {1/2}}$ and on the constant $\Lambda (\Omega)$  introduced in $\eqref{f:Lambda}$. Indeed, starting from the equality \eqref{f:charh}, some straightforward computations give
 $$ h_2 (\Omega)  |\Omega| ^ {3/2} = \Phi (I (\Omega), \Lambda (\Omega))\, , $$ 
 where
 $$\Phi (x, y):=\frac{x-\frac{2 (y-\pi
   )}{\sqrt{x^2-3 y+3 \pi
   }+x}}{\left(1-\frac{y-\pi
   }{\left(\sqrt{x^2-3 y+3 \pi
   }+x\right)^2}\right)^2}\,.$$
 We observe that $\Lambda (\Omega)$ and $I (\Omega)$ obey the inequalities
$$I ^ 2 (\Omega) \geq 4 \Lambda (\Omega) \geq 4 \Lambda (\Omega ^*_n) = 4 n \tan \big ( \frac{\pi}{n} \big ) ( > 4 \pi) \,. $$
The former is the isoperimetric inequality for convex polygons (see for instance \cite{CFG02}), and becomes an equality when the polygon is circumscribed to a disk; 
the latter, recalling the definition of $\Lambda (\Omega)$, comes  from the convexity of the map $t \mapsto \cot (t)$ on the interval $(0, \pi/2)$, and becomes an equality whan all the inner angles of the polygon are equal. 

We are thus led to minimize the function $\Phi (x, y)$ on the admissible region
$$\mathcal A_n :=  \Big \{ (x, y) \ :\  x^ 2 \geq 4 y \geq 4 n \tan \big (\frac{\pi}{n} \big )\Big \} \,.$$

We claim that the map $y \mapsto \Phi (x, y)$ is decreasing. Indeed, let us show that
\begin{equation}\label{f:deriy}
\frac{\partial \Phi (x, y)} {\partial y} = -\frac{x^3+(x^2 +12 \pi - 12 y)\sqrt{x^2-3 y+3\pi }}{4
   \left(x^2-4 y+4 \pi
   \right)^2} \leq 0  \qquad \text{ if } (x, y) \in \mathcal A _n \,,
   \end{equation}
or equivalently that 
\begin{equation}\label{f:psi} 
\Psi (x, y):= x^3+(x^2 +12 \pi - 12 y)\sqrt{x^2-3 y+3\pi } \geq 0  \qquad \text{ if } (x, y) \in \mathcal A _n \,.
\end{equation}
Since
$$ \frac{\partial \Psi (x, y)}{\partial y } = -\frac{27 \left(x^2-4 y+4 \pi
   \right)}{2 \sqrt{x^2-3 y+3
   \pi }} \leq 0 \qquad \text{ if } (x, y) \in \mathcal A  _n \,,$$
in order to obtain \eqref{f:psi} it is enough to show that
$$\eta (x):= \Psi  \big (x, \frac{x^2}{4} \big ) = x^3+(6
   \pi-  x^2)  \sqrt{x^2+12 \pi } \geq 0 \qquad \text{ if } x \geq 2 \sqrt{n \tan \big (\frac{\pi}{n} \big )} \,.$$
The latter inequality is readily checked, since the function $\eta$ turns out to be monotone decreasing, 
and satisfies the following asymptotic expansion as $x \to + \infty$:
$$\frac{\eta(x)}{x^3} \sim  \frac{54 \pi ^ 2 } {x^4} + o \Big ( \frac{1}{x^4} \Big )\,. $$ 

We have thus proved \eqref{f:deriy}, yielding 
$$\Phi (x, y) \geq \Phi \big (x, \frac{x ^ 2}{4} \big )= 
\frac{\left(\sqrt{x^2+12 \pi }+2 x\right)^3}{16 \left(x
   \left(\sqrt{x^2+12 \pi }+x\right)+4 \pi \right)} =:\zeta (x)
 \,.$$
Next we observe that the map $x \mapsto \zeta (x)$ is increasing  for $ x \geq 2 \sqrt {n \tan \big (\frac{\pi}{n} \big )} $. 
Indeed, we have
$$\zeta' (x) = \frac{3 \left(x
   \left(\sqrt{x^2+12 \pi
   }-x\right)+12 \pi
   \right)}{64 \pi } \geq 0 \qquad \text{ if } 
x \geq 2 \sqrt \pi\,.$$
We conclude that the minimum of $\Phi (x, y)$ over the region $\mathcal A_n$ is attained at 
$$(x_n, y _n) :=\Big ( 2 \sqrt {n \tan \big (\frac{\pi}{n} \big )}, n \tan \big (\frac{\pi}{n}) \Big )\,,$$
corresponding to the case when the convex polygon $\Omega$ is at the same time circumscribed to a disk and with all the inner angles equal, 
that is, $\Omega$ is the regular polygon. 
Accordingly, the expression of $\gamma (n)$ in \eqref{f:gamma} is found by evaluating $\Phi$ at $(x_n, y _n)$. 

\qed

\bigskip
Finally, we arrive at the following honeycomb-type result, which extends to the case of the $2$-Cheeger constant Corollary 9  in \cite{bfvv17}:

\begin{proposition}\label{t:h2}
There holds
\begin{equation}\label{parti1}
\lim _{k \to + \infty} \frac{|\Omega| ^{3/2}}{k ^ {5/2}} \inf  \Big \{  \sum_{i = 1} ^k  h _2 (E _i)  \ :\  \{E_i\} \in \C _k (\Omega) \Big \}  = h _ 2 (H) \, , 
\end{equation}  
where $h_2 ( H)$ denotes the $2$-Cheeger constant of the unit area regular hexagon. 
\end{proposition}

\proof 
The equality \eqref{parti1} 
follows by applying Theorem 2 
in \cite{bfvv17}. One has just  check that assumptions (H1), (H2), (H3) therein are fulfilled. 
Assumptions (H1) and (H2) are satisfied, because the map $\Omega \mapsto h _ 2 (\Omega)$ is monotone decreasing 
under inclusions and homogeneous of degree $-3$ under domain dilations. It remains to check assumptions (H3): according to Remark 4 (ii) in \cite{bfvv17}, 
in view of Proposition \ref{p:discreteFK}, it is enough to check  that the map $n \mapsto \gamma (n) ^ {2/5}$ 
admits a decreasing and convex extension on $[2, + \infty)$.  This can be done by elementary computations by exploiting the explicit expression of $\gamma (n)$ given in \eqref{f:gamma}. 
\qed

 \begin{corollary}\label{p:Cheeger}
There holds
\begin{equation}\label{Cheegersum2}
\lim _{k\to + \infty} \frac{|\Omega|^ {3/2} }{ k ^ {5/2} }   \inf  \Big \{  \sum_{i = 1} ^k   \frac{|\partial E _i| }{|E _i |^2}  \ :\  \{E_i\} \in \C _k (\Omega) \Big \} =  h _2( H) \, , 
\end{equation}
where $h_2 ( H)$ denotes the $2$-Cheeger constant of the unit area regular hexagon. 
\end{corollary}

\proof Set $m _k (\Omega)$ and $\widetilde m _k (\Omega)$ the infima at the r.h.s. of  \eqref{Cheegersum2} and \eqref{parti1} respectively. 
The corollary follows straigthforward from Proposition \ref{t:h2} by noticing that  $m _k (\Omega)=\widetilde m _k (\Omega)$.
Indeed, let $\{E_i \} \in \mathcal C _k (\Omega)$. The inequality $h _2(E _i) \leq |\partial E _i| / |E _i| ^2$ yields 
immediately $\widetilde m _k (\Omega) \leq m _k (\Omega)$.
Conversely, since $h_2 ( E_i ) = |\partial C (E_i) | / |C(E_i)|^2$, there holds $\sum _i h_2 ( E_i) \geq m _k (\Omega)$,
which yields 
$\widetilde m _k (\Omega) \geq  m _k (\Omega)$.
\qed

\section{Some intermediate results}\label{sec:prel}

We proceed separately in the cases of $\lambda _1 (\Omega, \beta)$ and $\tau ^ {-1} (\Omega, \beta)$. The results of this section heavily rely on some works by Giorgi-Smits and Sperb in the former case, and by Keady-McNabb in the latter case (see \cite{Giorgi, Giorgi2, Keady,Sperb}).

\subsection{Preliminaries to the proof of Theorem \ref{t:honeycomb}}

\begin{remark} It will be useful to keep in mind the following scaling law, 
which can be easily checked by change of  variables:
\begin{equation}\label{scaling}\l ( t E, \beta ) = \frac{1}{t ^ 2} \l ( E, t \beta ) \qquad \forall t >0\,.
\end{equation}
\end{remark}

\begin{proposition}\label{p:bounds}[upper and lower bounds for $\l (\Omega, \beta)$] Let $E$ be an open bounded Lipschitz set, and let
$\mu _2 (E)$ denote the first nonzero eigenvalue of the Neumann Laplacian in $E$. 
There holds:
\begin{eqnarray}
&  \displaystyle \l (E, \beta) \leq \beta \frac{|\partial E| } {|E |} \qquad \forall \beta \in \R \setminus \{0 \}
 & \label{upperbo}
\\ \noalign{\bigskip}
 &\displaystyle  \l (E, \beta) \geq \frac{1}{\frac{1}{\mu _ 2 (E) } + \frac{ |E|}{ \beta |\partial E |}}
 \qquad \forall \beta >0
 \label{lowerbo}
\\ \noalign{\bigskip}
&\displaystyle  \lim _{\beta \to 0} \frac{\l (E, \beta)}{\beta} =  \frac{|\partial E| } {|E |} \,.  \hskip 2.5 cm 
& \label{asymptotic}
\end{eqnarray}

\end{proposition}
\proof The upper bound \eqref{upperbo} is trivially obtained by taking as a test function $u \equiv 1$ in the definition of $\l(E, \beta)$. 
The lower bound \eqref{lowerbo} is due to Sperb, see \cite{Sperb}. The asymptotic behaviour \eqref{asymptotic} is a direct consquence of \eqref{upperbo}-\eqref{lowerbo} in case $\beta \rightarrow 0_+$ but  requires a further  control from below in case $\beta \rightarrow 0_-$, see for instance \cite[eq.(5)]{Giorgi}. \qed

\begin{proposition}\label{width}[estimate of $\l (\Omega, \beta)$ by the width] Assume $\beta >0$. Let $E$ be an open bounded convex set, and let $w ( E, \xi )$ denote the width of $E$ in some fixed direction $\xi$ ({\it i.e.}, the distance between two support planes of $E$ orthogonal to $\xi$). For every $\d>0$, there exists a positive constant $C_1 = C _1 (\beta, \delta)$ such that
$$w (E, \xi) < \delta \quad \Rightarrow \quad \l (E, \beta)  \geq\frac{C_1}{w(E, \xi)}\,.$$
\end{proposition}

\proof Throughout the proof, we write for brevity $w$ in place of $w ( E, \xi)$. We proceed in two steps.  First we obtain the inequality 
\begin{equation}\label{lowerbound} 
\l (E, \beta) \geq \frac{1}{w ^ 2} \l ( I, w \beta)\end{equation}
where $I = ( - 1/2, 1/2)$ denotes the unit interval of the real line (and its Robin eigenvalue is meant in dimension $1$) 
and then we show that the quotient
\begin{equation}\label{control}
\frac{\l ( I, w \beta)}{w}
\end{equation} 
admits a  positive finite  limit as $w \to 0 ^+$. 

To prove 
\eqref{lowerbound}, we slice $E$ in the direction $\xi$. Namely, we denote by $E _\xi$ the projection of $E$ onto the direction $\xi ^ \perp$, 
and for every $x \in E_\xi$ we set $(a_x, b _x):= E \cap (x + \R \xi)$. By Fubini's Theorem, if $u$ is a first eigenfunction for $\lambda _ 1 (E, \beta)$, we have
$$\begin{array}{ll} \lambda _ 1 (E, \beta)  & \displaystyle = \frac{ \int_{E} |\nabla u | ^ 2 + \beta \int_{\partial E} u ^ 2 }{ \int _{E} u ^ 2}
\\ \noalign{\bigskip}
& \displaystyle \geq  \frac{   \int _{E _{\xi}}  \int_{a_x} ^ {b_x}  |\nabla u | ^ 2 (x, y) \, dy \,  \, dx + \beta \int_{E_\xi} [u ^ 2(a_x) + u ^ 2 (b _x) ]  \, dx }{ \int _{E _{\xi}}  \int_{a_x} ^ {b_x}   u ^ 2(x, y) \, dy \, dx } 
\\ \noalign{\bigskip}
& \displaystyle \geq \min _{x \in E _{\xi}} \frac{ \int_{a_x} ^ {b_x}  |\nabla u | ^ 2 (x, y) \, dy + \beta  [u ^ 2(a_x) + u ^ 2 (b _x) ] } {  \int_{a_x} ^ {b_x}   u ^ 2(x, y) \, dy  } 
\\ \noalign{\bigskip}
& \displaystyle \geq \min _{x \in E _{\xi}} \l ( (a_x, b _x) , \beta) 
\geq  \l  \big ( \big ( - \frac{w}{2} , \frac{w}{2}  \big ) , \beta \big ) 
 = \frac{1}{w ^ 2}  \l ( I , w \beta) \,.
\end{array}                   
$$
where in the last line we have used the decreasing monotonicity of the map $B \mapsto \lambda _ 1(B, \beta)$ holding for balls $B$ for any $\beta>0$, and the scaling property \eqref{scaling}. 

To compute the limit of the quotient in \eqref{control},  we solve the b.v.p. which defines $ \lambda := \l (I, \alpha)$, that is we search for an even function on $I$ which satisfies
$$
\begin{cases}
-  u''  = \lambda  u  \quad \text{ in }  I &
\\
u' \big ( \frac{1}{2} \big ) + \alpha u   \big ( \frac{1}{2} \big ) = 0\,.  &
\end{cases}
$$
We have $u (x) = \cos (\sqrt \lambda x)$, and imposing the boundary condition we get  the following relation between $\alpha$ and $\lambda$:
\begin{equation}\label{relation}
\frac{\alpha ^ 2}{\lambda + \alpha ^ 2 } = \sin  ^2 \Big ( \frac{\sqrt \lambda }{2} \Big ) \,.
\end{equation}
In the limit as $\alpha \to 0$, we have $\lambda \to 0$; moreover, 
dividing \eqref{relation} by $\lambda$ and passing to the limit as $\alpha \to 0$, we get 
\begin{equation}\label{relation2}
\lim _{\alpha \to 0} \frac{\alpha ^ 2}{\lambda (\lambda + \alpha ^ 2 )} =  \frac{\sin  ^2 \Big ( \frac{\sqrt \lambda }{2} \Big )}{\lambda} = \frac{1}{4} \,.
\end{equation}
We observe that  $\alpha ^ 2 = o (\lambda)$; indeed, using \eqref{relation}, we see that
$$\frac{4 \alpha ^ 2}{\lambda} \sim  \frac{\alpha ^ 2} {\sin  ^2 \Big ( \frac{\sqrt \lambda } {2}\Big ) } = {\lambda + \alpha ^ 2} \to 0 \,.$$ 
 Therefore, \eqref{relation2} can be rewritten as
$$
\lim _{\alpha \to 0} \frac{\alpha ^ 2}{\lambda ^ 2 + o (\lambda ^ 2) } =  \frac{1}{4} \,,
$$ 
yielding
$
\lim _{\alpha \to 0} \frac{\alpha }{\lambda } = \frac{1}{2} 
$. 
We conclude that
$$\lim _{w \to 0} \frac{\l ( I, w \beta)}{w} =  2 \beta.
$$ 
\qed

\subsection{Preliminaries to the proof of Theorem \ref{t:honeycomb2}} \label{sec:prel2}

\begin{remark} The Robin torsion satisfies a scaling law analogue to \eqref{scaling}, which in this case reads:  
\begin{equation}\label{scaling2}\tau ^ {-1} ( t E, \beta ) = \frac{1}{t ^ 4} \tau ^ {-1} ( E, t \beta ) \qquad \forall t >0\,.
\end{equation}
\end{remark}

\begin{proposition}\label{p:bounds2}[upper and lower bounds for $\tau ^ {-1} (\Omega, \beta)$] Let $E$ be an open bounded Lipschits set with unit outer normal $\nu$, and let
$\Sigma _\infty (E) := \int _ E u _\infty$, being $u _\infty$ the unique solution  to the boundary value problem
$$\begin{cases}
- \Delta u = 1 & \text{ in } E
\\ 
u_ { \nu} = -  \frac{| E|}{|\partial E|}   & \text{ on } \partial E
\\
\int_{\partial E} u = 0 \,.
\end{cases}
$$
There holds:
\begin{eqnarray}
&  \displaystyle \tau ^ {-1} (E, \beta) \leq \beta \frac{|\partial E| } {|E |^2} \qquad \forall \beta \in \R \setminus \{0 \}
 & \label{upperbo2}
\\ \noalign{\bigskip}
 &\displaystyle   \tau ^ {-1} (E, \beta) \geq \frac{1}{\Sigma _\infty (E)  + \frac{ |E|^2}{ \beta |\partial E |}}
 \qquad \forall \beta >0
 \label{lowerbo2}
\\ \noalign{\bigskip}
&\displaystyle  \lim _{\beta \to 0} \frac{\tau ^ {-1}(E, \beta)}{\beta} =  \frac{|\partial E| } {|E |^2} \,.  \hskip 2.5 cm 
& \label{asymptotic2}
\end{eqnarray}

\end{proposition}

\bigskip
For the proof of Proposition \ref{p:bounds2} we need the following result;  similar statements for the first Robin eigenvalue can be found in \cite[Lemma 1]{AB} and \cite[Lemma 2.2]{Giorgi2}. 

\begin{lemma}\label{l:GS}
Let $E$ be an open bounded Lipschitz domain. 
\begin{itemize}
\item[(i)] There exists a constant $M _E$ such that
$$
\int _{\partial E} u ^ 2 \leq \int _E |\nabla u| ^ 2 + M _E \Big ( \int _E |u|  \Big ) ^ 2 \qquad \forall u \in H ^ 1 (E) \,.
$$
\item[(ii)] One can find $\eta_{E} >0$ such that, for every $\eta \geq \eta _{E}$, there  exists a constant $C (\eta) >0$, infinitesimal as $\eta \to + \infty$,  such that
$$\int_{\partial E} u ^ 2  \leq \eta \int_E  |\nabla u| ^ 2 + \frac{|\partial E|}{|E|^2} \big ( 1 + C(\eta) \big ) \Big (\int _{E} |u| \Big ) ^2   \qquad \forall u \in H ^ 1 (E)\,.$$
\end{itemize}
\end{lemma}

\proof
(i) By Lemma 1 in \cite{AB}, there exists  a constant $C>0$ such that
$$
\int _{\partial E} u ^ 2 \leq  \frac{1}{2} \int _E |\nabla u| ^ 2 +  C \int _E u^2   \qquad \forall u \in H ^ 1 (E) \,.
$$
Then it is enough to show the following claim: for any given $C>0$, there exists $M>0$ sufficiently large such that 
\begin{equation}\label{f:claimR}
 C \int _E u^2 \leq  \frac{1}{2} \int _E |\nabla u| ^ 2 + M  \Big ( \int _E |u|  \Big ) ^ 2 \qquad \forall u \in H ^ 1 (E) \,.
 \end{equation}
The claim is readily cheked by contradiction. Assume there exists a sequence $M _n \to + \infty$ such that, for every $n \in \N$, there exists $u _n \in H ^ 1 (E)$, with $\int _ E u _n ^ 2 = 1$, satisfying 
\begin{equation}\label{contr}
C  \geq  \frac{1}{2} \int _E |\nabla u_n| ^ 2 + M_n  \Big ( \int _E |u_n|  \Big ) ^ 2 \,.
\end{equation}
 Then the sequence $\{ u _n \}$ turns out to be bounded in $H ^ 1 (E)$ so that, up to passing to a (not relabeled) subsequence, it converges to some function $u$ weaky in $H ^ 1 (E)$ and strongly both in $L ^ 1 (E)$ and $L ^ 2 (E)$. Recalling that $\int _ E u _n ^ 2 = 1$ for every $n$, we find $\int _E u ^ 2 = 1$, whereas 
recalling that $M _n \to + \infty$, the inequality \eqref{contr} implies $\int _E |u| = 0$, contradiction. 

\smallskip
(ii) We proceed by contradiction. If statement (ii) is false, we can find $\delta >0$ and a sequence $\{u _n \} \subset H ^ 1 (E)$ such that, for every $n \in \N$, 
$$\int_{\partial E} u _n ^ 2  \geq n \int_E  |\nabla u_n| ^ 2 + \frac{|\partial E|}{|E|^2} ( 1 + \delta  ) \Big (\int _{E} |u_n| \Big ) ^2\,.$$
 The above inequality implies in particular that $ \int_E  |\nabla u_n| ^ 2 \neq 0$ for every $n \in \N$ (otherwise $u _n$ is constant and we get a contradiction since $\delta >0$). 
Then we can define the functions $v _n := \frac{u _n}{\big (\int _E|\nabla u _n | ^ 2 \big ) ^ {1/2}}$, which  satisfy
\begin{equation}\label{f:uno}
\int_{\partial E} v _n ^ 2  \geq n  + \frac{|\partial E|}{|E|^2} ( 1 + \delta  ) \Big (\int _{E} |v_n| \Big ) ^2\,.
\end{equation}
On the other hand, from statement (i), and taking into account that $\int _E |\nabla v_n| ^ 2 = 1$,  we know that 
\begin{equation}\label{f:due}
\int _{\partial E} v_n ^ 2 \leq 1 + M _E \Big ( \int _E |v_n|  \Big ) ^ 2 \,.
\end{equation}
Combining the two inequalities \eqref{f:uno} and \eqref{f:due}, we obtain a contradiction concluding the proof,  
provided we are able to show that the sequence $\int _E |v_n|$ remains bounded. 
Assume this is not the case. Then,  the sequence $w _n:=  \frac{v _n}{\int _E| v _n |  }$  satisfies $\int _E |\nabla w _n | ^ 2 \to 0$. 
Since $\int _E |w_n | = 1$, by exploiting claim \eqref{f:claimR} obtained above in the proof of statement (i), we infer that that $\{w_n \} $ remains bounded in $H ^ 1 (E)$. Hence, up to subsequences, it converges weakly in $H ^ 1 (E)$ and strongly in $L ^ 1 (E)$ to some limit $w$ which satisfies $\int _E |w| = 1$ and $\int _E |\nabla w|^2 = 0$. It follows that $w$ is uniquely determined as $\frac{1}{{|E|}}$. By uniqueness of this limit, the whole sequence turns out converge strongly in $H ^ 1 (E)$ to $w$, and hence also strongly in $L ^ 2 (\partial E)$. Then \eqref{f:uno} implies 
$$\frac{|\partial E|}{|E|^2} = \lim _n  \int_{\partial E} w _n ^ 2  \geq   \frac{|\partial E|}{|E|^2} ( 1 + \delta  ) \, , 
$$
contradiction.

\qed

\bigskip

{\it Proof of Proposition \ref{p:bounds2}}.   For the upper bound \eqref{upperbo2}, it's enough to take $u \equiv 1$ in the definition of $\tau ^ {-1}(E, \beta)$. 
The lower bound \eqref{lowerbo2} is due to Keady-McNabb, see \cite[inequality (4.9)]{Keady}. 
Let us prove the asymptotic behaviour \eqref{asymptotic2}. The limit as $\beta \to 0 ^ +$ is obtained immediately by combining the bounds \eqref{upperbo2} and \eqref{lowerbo2}. 
It remains to compute the limit as $\beta \to 0 ^ -$. 
The inequality $\liminf\limits _{\beta \to 0 ^ -} \frac{\tau ^ {-1} (E, \beta)  }{\beta} \geq \frac{|\partial E|}{|E|^2}$ follows immediately from  \eqref{upperbo2}. Let us show the converse inequality for the limsup. 
To that aim, we apply Lemma \ref{l:GS} by choosing $\eta = - \frac{1}{\beta}$ and taking as a function $u$ the solution to the Robin torsion problem, normalized so that $\int_E|u| = 1$. We obtain
$$\frac{\tau ^ {-1} (E, \beta)  }{\beta} =\int_{\partial E} u ^ 2 + \frac{1}{\beta}  \int_{E} |\nabla u| ^ 2  \leq   \frac{|\partial E|}{|E|^2} 
\big ( 1 + C \big (- \frac{1}{\beta} \big ) \big )  \Big (\int _{E} |u| \Big ) ^2 =    \frac{|\partial E|}{|E|^2} 
\big ( 1 + C \big (- \frac{1}{\beta} \big ) \big )   \,.$$
It follows that
$\limsup\limits _{\beta \to 0 ^ -} \frac{\tau ^ {-1} (E, \beta)  }{\beta} \leq \frac{|\partial E|}{|E|^2}$ as required. 
 \qed

\begin{proposition}\label{width2}[estimate  by the width] Assume $\beta >0$. Let $E$ be an open bounded convex set, and let $w ( E, \xi )$ denote the width of $E$ in some fixed direction $\xi$ ({\it i.e.}, the distance between two support planes of $E$ orthogonal to $\xi$). For every $\d>0$, there exists a positive constant $C_1 = C _1 (\beta, \delta)$ such that
$$w (E, \xi) < \delta \quad \Rightarrow \quad \tau ^ {-1}(E, \beta)  \geq\frac{C_1}{w^3(E, \xi)}\,.$$
\end{proposition}

\proof We proceed in a similar way as in the proof of Proposition \ref{width}. We still set $w:=w ( E, \xi)$, and $I := ( - 1/2, 1/2)$, and we   proceed again in two steps, showing first that  
\begin{equation}\label{lowerbound2} 
\tau ^ {-1} (E, \beta) \geq \frac{1}{w ^ 4 |E_\xi|} \tau ^ {-1}  ( I, w \beta)\end{equation}
and second that  the quotient
\begin{equation}\label{control2}
\frac{\tau ^ {-1}( I, w \beta)}{w}
\end{equation} 
admits a  positive finite  limit as $w \to 0 ^+$. 

To prove \eqref{lowerbound2}, we proceed by slicing. With the same notation as in the proof of Proposition \ref{width}, by using Fubini's Theorem and H\"older inequality, if $u$ is the solution to the Robin torsion problem, we have
$$\begin{array}{ll} \tau ^ {-1} (E, \beta)  & \displaystyle = \frac{ \int_{E} |\nabla u | ^ 2 + \beta \int_{\partial E} u ^ 2 }{ \Big ( \int _{E} |u| \Big ) ^ 2}
\\ \noalign{\bigskip}
& \displaystyle \geq  \frac{   \int _{E _{\xi}}  \int_{a_x} ^ {b_x}  |\nabla u | ^ 2 (x, y) \, dy \,  \, dx + \beta \int_{E_\xi} [u ^ 2(a_x) + u ^ 2 (b _x) ]  \, dx }{\Big (  \int _{E _{\xi}}  \int_{a_x} ^ {b_x}   |u (x, y)| \, dy \, dx \Big ) ^2} 
\\ \noalign{\bigskip}
& \displaystyle \geq \frac{1}{|E_\xi|} \min _{x \in E _{\xi}} \frac{ \int_{a_x} ^ {b_x}  |\nabla u | ^ 2 (x, y) \, dy + \beta  [u ^ 2(a_x) + u ^ 2 (b _x) ] } {\Big (  \int_{a_x} ^ {b_x}   |u (x, y)| \, dy   \Big ) ^ 2} 
\\ \noalign{\bigskip}
& \displaystyle \geq \frac{1}{|E_\xi|}  \min _{x \in E _{\xi}} \tau ^ {-1} ( (a_x, b _x) , \beta) 
\geq  \frac{1}{|E_\xi|}   \tau ^ {-1}  \big ( \big ( - \frac{w}{2} , \frac{w}{2}  \big ) , \beta \big ) 
 = \frac{1}{w ^ 4 |E_\xi|}  \l ( I , w \beta) \,.
\end{array}                   
$$
where in the last line we have used the decreasing monotonicity of the map $B \mapsto  \tau ^ {-1}(B, \beta)$ holding for balls $B$ for any $\beta>0$, and the scaling property
\eqref{scaling2}. 

To compute the limit of the quotient in \eqref{control2},  we solve the b.v.p. which defines $ \tau ^ {-1} (I, \alpha)$, that is we search for an even function on $I$ which satisfies
$$
\begin{cases}
-  u''  =  1  \quad \text{ in }  I &
\\
u' \big ( \frac{1}{2} \big ) + \alpha u   \big ( \frac{1}{2} \big ) = 0\,.  &
\end{cases}
$$
We find $u (x) = - \frac{x^2}{2} + \frac{1}{2 \alpha} + \frac{1}{8}$, with $\int _I u (x) \, dx = \frac{1}{12} + \frac{1}{2\alpha}$. 

We conclude that
$$\lim _{w \to 0} \frac{\tau ^ {-1} ( I, w \beta)}{w} = \lim _{w \to 0} \frac{1}{w}   \frac{1}{ \frac{1}{12} + \frac{1}{2 w \beta}} = 2 \beta \,.$$ 

\qed

\section {Proof of Theorem \ref{t:honeycomb}}\label{sec:proof}
 We start from the following fact: setting
\begin{equation}\label{Cheegersum}
m _{k} (\Omega) = \inf  \Big \{  \sum_{i = 1} ^k   \frac{|\partial E _i| }{|E _i |}  \ :\  \{E_i\} \in \C _k (\Omega) \Big \} 
\end{equation}
there holds
\begin{equation}\label{aCheegersum}
 \lim _{k\to + \infty} \frac{|\Omega|^ {1/2} }{ k ^ {3/2} } m _k (\Omega) =  h ( H) 
 \end{equation}
This is obtained immediately  by applying Corollary 9 
in \cite{bfvv17}, and arguing as in the proof of Corollary \ref{p:Cheeger}. 

We now proceed separately in the two cases $\beta >0$ and $\beta <0$, assuming 
without loss of generality that $|\Omega| = 1$.  
\bigskip

$\bullet$ {\bf Case $\beta >0$.}

\smallskip
In view of \eqref{upperbo} in Proposition \ref{p:bounds} and \eqref{aCheegersum}
, we have 

\begin{equation}\label{ub1} 
\limsup _{k\to + \infty} \frac{ r _k (\Omega, \beta) }{ k ^ {3/2} }  \leq
\limsup _{k\to + \infty} \beta \frac{ m _k (\Omega) }{ k ^ {3/2} } 
= \beta h ( H) 
\end{equation}

We are going to show that
\begin{equation}\label{tesi1} 
\liminf _{k\to + \infty} \frac{ r _k (\Omega, \beta) }{ k ^ {3/2} }  \geq
\liminf _{k\to + \infty} \beta \frac{ m _k (\Omega) }{ k ^ {3/2} } 
= \beta h ( H)
\end{equation}

By \eqref{ub1} 
, we can choose $\overline k$ sufficiently large so that, for every $k \geq \overline k$, it holds
\begin{equation}\label{stima}
\frac{ r _k (\Omega, \beta) }{ k ^ {3/2} } \leq 2 \beta h ( H) \,.
\end{equation}
Then, for every $k \geq \overline k$, we let $\{ \omega^*_1, \dots, \omega ^* _k \}$ 
be a convex cluster in $\Omega$ such that
\begin{equation}\label{u1}
 \sum_{i = 1} ^ k  \l (\omega^* _i, \beta)  \leq  \displaystyle r _{k} (\Omega, \beta) +1 
 \end{equation}

For a given $\e>0$,  we introduce the class of convex bodies such that  the ratio between the width  in a direction orthogonal to a diameter and the diameter  is at least $\e$.  We denote it by 
$$\conv (\e) := \Big \{ E \in \K ^ 2 \ :\ \frac{w (E, \xi)}{{\rm diam}(E)} \geq \e \ \text{ for some   $\xi \in S ^1$ orthogonal to a diameter} \Big \}\,. $$

Then we consider the following families of indices associated with the clusters  $\{ \omega^*_1, \dots, \omega ^* _k \}$ 
$$
\theta _{k, \e} := \big \{ i \in \{ 1, \dots, k \} \ :\ \omega^ * _i \in \conv (\e) \big \} \,, \qquad \theta _ {k,\e} ^ c := \big \{1, \dots, k  \big \} \setminus \theta _\e\,.
$$

\smallskip
We can estimate $r_k (\Omega; \beta)$ 
from below as follows
\begin{equation}\label{aux1}
1 + r _k (\Omega, \beta)  \geq \sum _{i=1} ^k \lambda _1 ( \omega ^*_i, \beta) \geq 
 \sum _{i \in \theta_{k,\e}}  \lambda _1 ( \omega ^*_i, \beta)  \,.
 \end{equation}
 
We are thus led to introduce the auxiliary problems
\begin{equation} \label{pbaux1}
r _{k, \e} (\Omega, \beta) := \inf \Big \{ 
\sum _{i \in \theta _{k, \e}}\lambda _1 ( E _i, \beta) \ :\  \{E _i \} \in \mathcal C _{\sharp \theta _{k,\e}} (\Omega) \, ,\  E _i \in \conv (\e) \Big \} \,.
\end{equation}

In order to show \eqref{tesi1},
we are going to exploit the lower bound \eqref{aux1},
and to estimate from below $r _{k, \e} (\Omega, \beta)$ 
in terms of the corresponding auxiliary problems 
\begin{equation}\label{mpbaux1}
m _{k, \e} (\Omega) := \inf \Big \{ 
\sum _{i \in \theta _{k, \e}}  \frac{|\partial E _i|} {|E_i |}  \ :\  \{E _i \} \in \mathcal C _{\sharp \theta _{k,\e}} (\Omega) \, ,\  E _i \in \conv (\e) \Big \} \,.
\end{equation}
We divide the remaining of the proof in three steps. 

\medskip
{\it Step I: for $k$ large enough,  it holds
\begin{equation}\label{claim1}
\sharp \theta _{k, \e} \geq ( 1 - C \e ^ {1/3} ) k
\,, 
\end{equation}
where  $C$ stands for a positive constant, not depending on $k$ nor on $\e$.
Consequently, 
\begin{equation}\label{behave} 
\lim _{k \to  + \infty} \frac{m _{k, \e} (\Omega)}{ (\sharp \theta_{k, \e}) ^ {3/2}   } = h ( H)\,.
\end{equation}
}

\smallskip 
Let us first observe that \eqref{behave} is a straightforward consequence of \eqref{claim1}. Indeed, by definition it is clear that $m _{k, \e}(\Omega) \geq m _k (\Omega)$
On the other hand, provided $\e < 1/2$,  we have $H \in \conv (\e)$, so that a configuration which is asymptotically optimal, in the limit as $k \to + \infty$,  for $m _k (\Omega)$ 
is admissible for $m _{k, \e} (\Omega)$.
This yields \eqref{behave} since, by \eqref{claim1}, we know that $\sharp \theta _{k, \e}\sim k$ 
as  $k \to + \infty$.

In order to estimate the cardinality of $\theta  _{k,\e}$,
we first obtain a bound on the width of the cells $\omega ^ * _i$ and $\Omega^ *_i$. 
Hereafter we denote for brevity
$d _i$ and  $w _i$ the diameter of such cells, and their width in the direction orthogonal to a diameter. 

We have
\begin{equation} \label{prima1}
 1 = |\Omega| \geq \sum _{i = 1}^ k |\omega^*_i| \geq \frac{1}{2}  \sum _{i = 1}^ k  d _i w _i \geq \frac{1}{2} \sum _{i \in \theta ^ c _\e} \frac{w_i ^ 2}{\e} 
\end{equation}
where we have used the fact that, by convexity, the area of each cell is bounded from below by  $(d_i w _i)/2$ and the fact that, by definition, for cells in $\theta^ c _{k,\e}$, it holds  $d _i \geq w _i/\e$. 

Now we must proceed  to estimate the cardinality  of  $\theta  _{k,\e}$.

Starting from \eqref{prima1} and using Proposition \ref{width},  the elementary inequality between the $2$-mean and the $(-1)$-mean, inequality \eqref{u1}, and inequality \eqref{stima}, we obtain, for $k$ large enough:

$$\begin{array}{ll}
1 & \displaystyle \geq \frac{1}{2 \e} 
\sum _{i \in \theta _{k, \e} ^ c} w _i ^ 2 
\geq \frac{C_1^2}{2 \e  }   
 \sum _{i \in \theta _{k, \e} ^ c} \Big ( \frac{1}{\lambda _1 (\omega ^*_i, \beta) }\Big ) ^ 2 \
\\ \noalign{\bigskip} 
& \displaystyle
 \geq \frac{C_1^2}{2 \e  } \frac 
 {(\sharp \theta _{k, \e} ^ c) ^ 3 }  
{\displaystyle \Big (\sum _{i \in \theta _{k, \e} ^ c}  \lambda _1 (\omega ^*_i, \beta) \Big ) ^ 2 }
\geq \frac{C_1^2}{2 \e  } \frac 
 {(\sharp \theta _{k, \e} ^ c) ^ 3 }  
{\displaystyle \Big ( 3 \beta h (H) k ^ {3/2} \Big ) ^ 2 }\,.
\end{array}
$$
Hence 
inequality in \eqref{claim1} is satisfied, for $k$ large enough, with $C  =  \Big ( \frac{18 \beta ^ 2 h ^ 2 (H)}{C_1^2} \Big ) ^ {1/3}$.

\medskip
{\it Step II: For $k$ large enough,  the infima $r_{k, \e} (\Omega, \beta)$
and $m_{k, \e} (\Omega, \beta)$
introduced in \eqref{pbaux1} and \eqref{mpbaux1} satisfy
\begin{equation}\label{claim2a} 
\liminf_{k \to + \infty}\frac{r _{k, \e} (\Omega, \beta) }{k ^ {3/2}}    \geq (1- \e) \beta   \liminf_{k \to + \infty}  \frac{ m _{k, \e} (\Omega)}{k ^ {3/2}} \,.
\end{equation}}

Let $\{E _i \}$ be a cluster in $\mathcal C _{\sharp \theta _{k,\e}} (\Omega)$,
with $E _i  \in \conv (\e)$. 
In order to estimate $\lambda _ 1 ( E_i, \beta)$, we 
introduce the following constant depending only on $\e$:
$$K _\e := \inf \Big \{ \mu _2 (E) |E| \ :\ E \in \conv (\e) \Big \}\,.$$
Notice that $K _\e$ is strictly positive (and attained) because the family $\conv (\e)$ is closed in the Hausdorff topology. 
Then we distinguish the cells of the cluster $\{E _i \}$ into two disjoint subclasses, in which we are able to provide respectively a lower and an upper bound for $\lambda _1 ( E_i, \beta)$.

\medskip
\begin{itemize}
\item[Class 1]: cells with $\beta |\partial E _i|  \leq \e K _\e$. For such cells, it holds 
\begin{equation}\label{t1}
\lambda _ 1 ( E _i, \beta) \geq ( 1 - \e) \beta \frac{|\partial E _i| }{|E_i|}\,.
\end{equation}
\medskip
Namely,  using the lower bound \eqref{lowerbo} given by Proposition \ref{p:bounds}, we have
$$\lambda  _1 ( E _i,\beta) \geq  \frac{1}{\frac{1}{\mu _ 2 (E_i) } + \frac{ |E_i|}{ \beta |\partial E_i |}}
= \beta \frac{|\partial E _i| }{|E_i|} \Big ( 1 - \frac{ \beta |\partial E _i |}{\beta |\partial E _i | + \mu _2 ( E _i) |E_i|} \Big ) \,.
$$ 
Therefore, the required estimate \eqref{t1} is satisfied provided 
$$\frac{ \beta |\partial E _i |}{\beta |\partial E _i | + \mu _2 ( E _i) |E_i|}  \leq \e\, , $$
which holds for cells of Class 1, as the inequality $\beta |\partial E _i| \leq \e K _ \e \leq \e \mu _ 2 ( E_i) |E_i|$ is in force. 
\medskip
\item[Class 2]: cells with $\beta |\partial E _i|  > \e K _\e$. For such cells, it holds  
\begin{equation}\label{t2}
\lambda _ 1 ( E _i, \beta) \leq \beta \frac{|\partial E _i|}{|E_i|} \leq \frac{32 \, \beta ^ 2 }{ \e ^ 2  K _ \e }\,. 
\end{equation}
Namely, using the upper bound \eqref{upperbo} given by Proposition \ref{p:bounds}, 
and the elementary estimates $|\partial E_i| \leq 4 d _i$, $|E _i| \geq \frac{1}{2} d _i w _i $ 
(being   $d_i$ the diameter of $E_i$, and $w_i$ its width in a direction orthogonal to a diameter),
we obtain
$$ \lambda  _1 ( E _i,\beta)  \leq \beta \frac{|\partial E _i|}{|E_i|} \leq \beta \frac{4 d _i}{\frac{1}{2} d _i w _i}  = \frac{ 8 \beta}{w_i} 
$$  
Therefore, the required estimate \eqref{t2} is satisfied provided 
\begin{equation}\label{nec2}
w _i \geq \frac{ \e ^ 2 K _\e}{4 \beta}
\end{equation}
which holds for cells of Class 2. Indeed for such cells the inequality $4 d _i \beta \geq \e K _ \e$ is in force, which yields $d _i \geq \frac{\e K _\e} {4 \beta}$. In turn, the latter inequality implies  \eqref{nec2} since $E _i \in \conv (\e)$. 
\end{itemize}

\bigskip
Now we proceed to prove the estimate in \eqref{claim2a}. 
Let  $\{E _i \}$ be a cluster in $\mathcal C _{\sharp \theta _{k,\e}} (\Omega)$, with $E _i  \in \conv (\e)$. We set for brevity   
$$
\begin{array}{ll}
& \theta ^ {(1)} _{k, \e} (\{E_i\}) := \Big \{ i \in \theta _{k, \e} \ :\ E _i \text{ is of Class 1} \Big \} 
\\ 
\noalign{\bigskip} 
& \theta ^ {(2)} _{k, \e} (\{E_i\}) := \Big \{ i \in \theta _{k, \e} \ :\ E _i \text{ is of Class 2} \Big \} \,.
\end{array} 
$$

We start by noticing that, in the limit as $k \to + \infty$, the infimum which defines $m _{k, \e} (\Omega)$ has the same asymptotic behaviour if we restrict the  the sum  to indices in the family  $\theta ^ {(1)} _{k, \e}(\{E_i\})$. More precisely, we have:  
\begin{equation}\label{asyb}
 \begin{array}{ll} 
& \displaystyle  \liminf _{k \to + \infty} \frac{1}{k ^ {3/2}} \inf  \Big \{ \sum _{i \in \theta  _{k, \e} }   \frac{|\partial E _i|}{ |E_i|}  \ :\    \{E _i \} \in \mathcal C _{\sharp \theta _{k,\e}} (\Omega) \, , \ E _i  \in  
\conv (\e) \Big \}  = 
\\ \noalign{\medskip} 
 & \displaystyle \liminf _{k \to + \infty} \frac{1}{k ^ {3/2}}  \inf  \Big \{ \sum _{i \in \theta ^ {(1)} _{k, \e}(\{E_i\}) }   \frac{|\partial E _i|}{ |E_i|}  \ :\    \{E _i \} \in \mathcal C _{\sharp \theta _{k,\e}} (\Omega) \, , \ E _i  \in  
\conv (\e) \Big \} \,.
\end{array} 
\end{equation}

Indeed, on one hand we know from \eqref{claim1} and \eqref{behave}  that
\begin{equation}\label{asy1}
 m _{k, \e} (\Omega) = \inf  \Big \{ \sum _{i \in \theta  _{k, \e} }   \frac{|\partial E _i|}{ |E_i|}  \ :\    \{E _i \} \in \mathcal C _{\sharp \theta _{k,\e}} (\Omega) \, , \ E _i  \in  
\conv (\e) \Big \}  \sim h ( H) k ^{3/2}\, ;
\end{equation}
on the other hand,  for any admissible cluster $\{E _i \}$ in $\mathcal C _{\sharp \theta _{k,\e}} (\Omega)$, with $E _i  \in \conv (\e)$, by \eqref{t2} we have
\begin{equation}\label{asy2} \sum _{i \in \theta ^ {(2)} _{k, \e}(\{E_i\}) }   \frac{|\partial E _i|} {|E_i|} \leq \frac{32 \beta }{\e ^ 2 K _\e} \, k\, , 
\end{equation} 
where the quantity $\frac{32 \, \beta  }{ \e ^ 2  K _ \e }$ is independent of $k$.  Then \eqref{asyb} follows by  \eqref{asy1} and \eqref{asy2}. 

Now, exploiting \eqref{t1} and \eqref{asyb}, we obtain 
$$\begin{array}{ll}
\displaystyle \liminf _{k \to + \infty} \frac{r _{k, \e} (\Omega, \beta) }{k ^ {3/2}}  & \displaystyle = \liminf _{k \to + \infty} \frac{1}{k ^ {3/2}}   \inf
\Big \{ \sum _{i \in \theta _{k, \e}  }  \lambda _1 ( E _i, \beta)   \ :\  \{E _i \} \in \mathcal C _{\sharp \theta _{k,\e}} (\Omega) \, , \ E _i  \in \conv (\e) \Big \} \\   \noalign{\bigskip} & \displaystyle \geq \liminf _{k \to + \infty} \frac{1}{k ^ {3/2}}  
\inf\Big \{  \sum _{ i \in \theta ^ {(1)} _{k, \e}(\{E_i\}) } \lambda _1 ( E _i, \beta)  \ :\  \{E _i \} \in \mathcal C _{\sharp \theta _{k,\e}} (\Omega) \, , \ E _i  \in \conv (\e) \Big \}
 \\   \noalign{\medskip} & \displaystyle \geq (1- \e) \beta   \liminf _{k \to + \infty} \frac{1}{k ^ {3/2}} 
\inf  \Big \{ \sum _{i \in \theta ^ {(1)} _{k, \e}(\{E_i\}) }  \frac{|\partial E _i|}{ |E_i|}  \ :\    \{E _i \} \in \mathcal C _{\sharp \theta _{k,\e}} (\Omega) \, , \ E _i  \in \conv (\e)\Big \}
\\   \noalign{\medskip} & \displaystyle = (1- \e) \beta   \liminf _{k \to + \infty} \frac{1}{k ^ {3/2}} 
\inf  \Big \{ \sum _{i \in \theta _{k, \e}}  \frac{|\partial E _i|}{ |E_i|}  \ :\    \{E _i \} \in \mathcal C _{\sharp \theta _{k,\e}} (\Omega) \, , \ E _i  \in \conv (\e)\Big \}
 \\   \noalign{\medskip} &= \displaystyle (1- \e) \beta  \liminf _{k \to + \infty} \frac{ m _{k, \e} (\Omega)}{k ^ {3/2}} . 
\end{array}
$$

\bigskip

{\it Step III: The lower bound \eqref{tesi1}
holds true}. 
 
 \bigskip 
By \eqref{aux1}, \eqref{claim1}, and \eqref{claim2a}, we have 
\begin{equation}\label{chain1}
\begin{array}{ll}  \displaystyle \liminf _{k \to + \infty}  \frac{1 + r _k (\Omega, \beta) }{k ^ {3/2}}  &\displaystyle =  \liminf _{k \to + \infty}  
\frac{r _{k, \e} (\Omega, \beta)}{k ^ {3/2} } \geq (1 - \e) \beta  \liminf _{k \to + \infty}  \frac{m _{k, \e} (\Omega)}{k ^ {3/2} }  
 \\ \noalign{\medskip} &
\displaystyle =  (1 - \e) \beta \liminf _{k \to + \infty}    \Big (  \frac{\sharp \theta _{k, \e}}{k} \Big )  ^ {3/2}  \frac{m _{k, \e} (\Omega)}{ (\sharp \theta _{k, \e}) ^ {3/2}   }   \\ \noalign{\medskip} &  \geq \displaystyle (1 - \e) \beta 
\big (  1 -  C \e ^ {1/3}  \big )  ^ {3/2}    \liminf _{k \to + \infty} \frac{m _{k, \e} (\Omega)}{ (\sharp \theta _{k, \e}) ^ {3/2}   }\,.  \end{array}
\end{equation}

By \eqref{chain1} and \eqref{behave}, we conclude that 
$$\liminf _{k \to + \infty}  \frac{ r _k (\Omega, \beta) }{k ^ {3/2}}  \geq  (1 - \e) \beta 
\big (  1 -  C \e ^ {1/3}  \big )  ^ {3/2} h ( H) \,.
$$
Eventually, in the limit as $\e \to 0 ^+$, we obtain \eqref{tesi1}. 

\bigskip
$\bullet$ {\bf Case $\beta <0$.} 

In view of inequality \eqref{upperbo} given by Proposition \ref{p:bounds} and the asymptotic equality \eqref{aCheegersum}, and taking into account that $\beta <0$, we have 

$$\begin{array} {ll} 
 \displaystyle \limsup _{k\to + \infty} \frac{r _k (\Omega, \beta) }{ k ^ {3/2} }  
& \displaystyle =  \limsup _{k\to + \infty} \frac{1 }{ k ^ {3/2} }   
\sup  \Big \{  \sum_{i = 1} ^ k  \l (E _i, \beta) \ :\  \{E_i\} \in \C _k (\Omega) \Big \}
\\ 
& \displaystyle \leq
\limsup _{k\to + \infty}\frac{1 }{ k ^ {3/2} }   
\sup  \Big \{  \sum_{i = 1} ^ k  \beta \frac{|\partial E _i|}{|E_i|} \ :\  \{E_i\} \in \C _k (\Omega) \Big \} \\ \noalign{\bigskip}
& \displaystyle =
\limsup _{k\to + \infty}\frac{\beta }{ k ^ {3/2} }   
\inf  \Big \{  \sum_{i = 1} ^ k  \frac{|\partial E _i|}{|E_i|} \ :\  \{E_i\} \in \C _k (\Omega) \Big \} \\ \noalign{\bigskip}
& \displaystyle =
\limsup _{k\to + \infty} \frac{ \beta \, m _k (\Omega) }{ k ^ {3/2} } = \beta h (H) \,.
\end{array}
$$
To prove the converse estimate,  since $|\Omega| = 1$, we observe that for every $\vps >0$ there exists $k_\vps$ such that for every $k\ge k_\vps$, the set $\Omega$ contains the convex $k$-cluster 
\begin{equation}\label{optimalcluster}
\big \{ \Big (\frac{1-\vps}{k}\Big )^\frac 12 {C(H)}  \, , \dots,  \Big (\frac{1-\vps}{k}\Big )^\frac 12{C(H)}   \big \}
\end{equation}
given by $k$ copies of the 
Cheeger set  $C (H)$ of the unit area regular hexagon, each one scaled so to have area $\frac{1-\vps}{k}$. Hence, using also the asymptotic behaviour \eqref{asymptotic} and the scaling property \eqref{scaling}, we get
$$\begin{array} {ll} 
 \displaystyle \liminf _{k\to + \infty} \frac{r _k (\Omega, \beta) }{ k ^ {3/2} }  
& \displaystyle =  \liminf _{k\to + \infty} \frac{1 }{ k ^ {3/2} }   
\sup  \Big \{  \sum_{i = 1} ^ k  \l (E _i, \beta) \ :\  \{E_i\} \in \C _k (\Omega) \Big \}
\\  \noalign{\smallskip}
& \displaystyle \geq
\liminf _{k\to + \infty}\frac{1 }{ k ^ {3/2} }   
k\, \lambda _ 1 \Big (  \Big (\frac{1-\vps}{k}\Big )^\frac 12 {C(H)} ,  \beta \Big )  \\ \noalign{\bigskip}
& \displaystyle
= \liminf _{k\to + \infty}\frac{1 }{ k ^ {3/2} }   \frac{k ^ 2  }{(1-\vps)}
\, \lambda _ 1 \Big (  C(H) , \frac{(1-\vps)^ {1/2} \beta}{k ^ {1/2}} \Big )   \\ \noalign{\bigskip}
& \displaystyle
= \liminf _{k\to + \infty}\frac{1 }{ k ^ {3/2} }   \frac{k ^ 2  }{(1-\vps)} 
  \frac{(1-\vps)^ {1/2} \beta}{k ^ {1/2}}   \frac {|\partial C (H)|} {|C (H)|} =  \frac{1 }{(1-\vps)^{1/2}}  \beta h (H)  \,. \\ \noalign{\bigskip}
\end{array}
$$
Since $\vps >0$ is arbitrary, the proof is concluded.
\qed 

\section {Proof of Theorem \ref{t:honeycomb2}}\label{sec:proof2}
 On the basis of the results established in Section \ref{sec:prel2}, the proof of Theorem \ref{t:honeycomb2} proceeds along the same line 
as the proof of Theorem \ref{t:honeycomb}. Hence we present it more concisely, often referring to the proof of Theorem \ref{t:honeycomb} whenever the two proofs 
are basically the same. We still assume that $|\Omega | = 1$, and we now set 

\begin{equation}\label{Cheegersum22}
m _{k} (\Omega) = \inf  \Big \{  \sum_{i = 1} ^k   \frac{|\partial E _i| }{|E _i |^2}  \ :\  \{E_i\} \in \C _k (\Omega) \Big \} 
\end{equation}

$\bullet$ {\bf Case $\beta >0$.}

\smallskip
In view of inequality \eqref{upperbo2} in Proposition \ref{p:bounds2} and Corollary \ref{p:Cheeger}, we have 

\begin{equation} \label{ub12}
\limsup _{k\to + \infty} \frac{ r _k (\Omega, \beta) }{ k ^ {5/2} }  \leq
\limsup _{k\to + \infty} \beta \frac{ m _k (\Omega) }{ k ^ {5/2} } 
= \beta h _2( H) 
\end{equation}

We have to prove that
\begin{equation}\label{tesi12}
\liminf _{k\to + \infty} \frac{ r _k (\Omega, \beta) }{ k ^ {5/2} }  \geq
\liminf _{k\to + \infty} \beta \frac{ m _k (\Omega) }{ k ^ {5/2} } 
= \beta h _2( H) \,.
\end{equation}

By \eqref{ub12},
we can choose $\overline k$ sufficiently large so that, for every $k \geq \overline k$, it holds
\begin{equation}\label{stima2}
\frac{ r _k (\Omega, \beta) }{ k ^ {5/2} } \leq 2 \beta h _2 ( H) \,.
\end{equation}
Then, for every $k \geq \overline k$, we let $\{ \omega^*_1, \dots, \omega ^* _k \}$ 
be a convex cluster in $\Omega$ such that
\begin{equation} \label{u12}
\sum_{i = 1} ^ k  \tau ^ {-1} (\omega^* _i, \beta)  \leq  \displaystyle r _{k} (\Omega, \beta) +1 \,.
\end{equation}

For a given $\e>0$,  we introduce the class of convex bodies $\conv (\e)$ defined as in the proof of Theorem \ref{t:honeycomb}, and accordingly we consider 
the families of indices $\theta _{k, \e}$
associated with the cluster  $\{ \omega^*_1, \dots, \omega ^* _k \}$ 
as done in such proof.

\smallskip
Also in the present setting, we have the lower bound
\begin{equation} \label{aux12}
1 + r _k (\Omega, \beta)  \geq \sum _{i=1} ^k \tau ^ {-1}  ( \omega ^*_i, \beta) \geq 
 \sum _{i \in \theta_{k,\e}}  \tau ^ {-1} ( \omega ^*_i, \beta) \,.
\end{equation}

Hence we introduce the auxiliary problems:
\begin{eqnarray}
& \displaystyle r _{k, \e} (\Omega, \beta) := \inf \Big \{ 
\sum _{i \in \theta _{k, \e}} \tau ^ {-1} ( E _i, \beta) \ :\  \{E _i \} \in \mathcal C _{\sharp \theta _{k,\e}} (\Omega) \, ,\  E _i \in \conv (\e) \Big \} & \label{pbaux12}
\\ \noalign{\medskip}
& \displaystyle m _{k, \e} (\Omega) := \inf \Big \{ 
\sum _{i \in \theta _{k, \e}}  \frac{|\partial E _i|} {|E_i |^2}  \ :\  \{E _i \} \in \mathcal C _{\sharp \theta _{k,\e}} (\Omega) \, ,\  E _i \in \conv (\e) \Big \}\,. & \label{mpbaux12}
\end{eqnarray} 
We divide the remaining of the proof in three steps. 

\medskip
{\it Step I: for $k$ large enough,  it holds
\begin{equation}\label{claim12}
\sharp \theta _{k, \e} \geq ( 1 - C \e ^ {3/5} ) k \, , 
\end{equation}
where  $C$ stands for a positive constant, not depending on $k$ nor on $\e$.
Consequently, 
\begin{equation}\label{behave2} 
\lim _{k \to  + \infty} \frac{m _{k, \e} (\Omega)}{ (\sharp \theta_{k, \e}) ^ {3/2}   } = h ( H)\,.
\end{equation}
}

\smallskip 
The equality \eqref{behave2} is deduced from \eqref{claim12} exactly in the same way as in the proof of Theorem \ref{t:honeycomb}.  
We proceed to prove the estimate \eqref{claim12}. Denoting by 
$d _i$ and  $w _i$ the diameter of the cells $\omega ^*_i$ or $\Omega ^*_i$, and their width in the direction orthogonal to a diameter,  we still have the inequality obtained in \eqref{prima1}, namely
\begin{equation}\label{prima1122} 
 1  \geq \frac{1}{2} \sum _{i \in \theta ^ c _\e} \frac{w_i ^ 2}{\e} 
\,.\end{equation}

\bigskip
{\it -- Estimate of $\sharp \theta _{k, \e}$}. 
Starting from 
\eqref{prima1122} and using Proposition \ref{width2},  the elementary inequality between the $\big ( \frac{2}{3}\big )$-mean and the $(-1)$-mean, inequality \eqref{stima2}, and inequality \eqref{u12}, we obtain, for $k$ large enough:

$$\begin{array}{ll}
1 & \displaystyle \geq \frac{1}{2 \e} 
\sum _{i \in \theta _{k, \e} ^ c} w _i ^ 2 
\geq \frac{C_1^{2/3}}{2 \e  }   
 \sum _{i \in \theta _{k, \e} ^ c}  \tau ^ {2/3} (\omega ^*_i, \beta)  \
\\ \noalign{\bigskip} 
& \displaystyle
 \geq \frac{C_1^{2/3}}{2 \e  } \frac 
 {(\sharp \theta _{k, \e} ^ c) ^ {5/3} }  
{\displaystyle \Big (\sum _{i \in \theta _{k, \e} ^ c}  \tau ^ {-1} (\omega ^*_i, \beta) \Big ) ^ {2/3} }
\geq \frac{C_1^{2/3}}{2 \e  } \frac 
 {(\sharp \theta _{k, \e} ^ c) ^ {5/3} }  
{\displaystyle \Big ( 3 \beta h _2(H) k ^ {5/2} \Big ) ^ {2/3} }\,.
\end{array}
$$
Hence
inequality \eqref{claim12} is satisfied, for $k$ large enough, with $C  =  \Big ( \frac{72 \beta ^ 2 h _2^ 2 (H)}{C_1^2} \Big ) ^ {1/5}$.

\medskip
{\it Step II: For $k$ large enough,  the infima $r_{k, \e} (\Omega, \beta)$,  
and $m_{k, \e} (\Omega, \beta)$, 
introduced in \eqref{pbaux12}-\eqref{mpbaux12}
satisfy
\begin{equation}\label{claim2a2}
\liminf_{k \to + \infty}\frac{r _{k, \e} (\Omega, \beta) }{k ^ {3/2}}    \geq (1- \e) \beta   \liminf_{k \to + \infty}  \frac{ m _{k, \e} (\Omega)}{k ^ {3/2}}  \, .
\end{equation}}

Let $\{E _i \}$ be a cluster in $\mathcal C _{\sharp \theta _{k,\e}} (\Omega)$,
with $E _i  \in \conv (\e)$. 
In order to estimate $\tau ^ {-1} ( E_i, \beta)$, we 
introduce the following constant depending only on $\e$:
$$K _\e := \inf \Big \{ (\Sigma _\infty  (E)  ) ^ {-1} |E|^2 \ :\ E \in \conv (\e) \Big \}\,,$$
the constant $\Sigma _\infty (E)$ being defined as in Proposition \ref{p:bounds2}.
We observe that $K _\e$ is strictly positive and attained because the family $\conv (\e)$ is closed in the Hausdorff topology, 
and  the functional $ (\Sigma _\infty  (E)  ) ^ {-1} |E|^2$ is continuous and scale-invariant (indeed, $\Sigma _\infty$ is easily seen from its definition to be homogeneous of degree $4$ under dilations). 
Then we distinguish the cells of the cluster $\{E _i \}$ into the same two subclasses as in the proof of Theorem \ref{t:honeycomb}.

\medskip
\begin{itemize}
\item[Class 1]: cells with $\beta |\partial E _i|  \leq \e K _\e$. For such cells, it holds 
\begin{equation}\label{t12}
\tau ^ {-1} ( E _i, \beta) \geq ( 1 - \e) \beta \frac{|\partial E _i| }{|E_i|^2}\,.
\end{equation}
\medskip
Namely,  using the lower bound \eqref{lowerbo2} given by Proposition \ref{p:bounds2}, we have
$$\tau ^ {-1} ( E _i,\beta) \geq  \frac{1}{\Sigma _\infty (E_i)  + \frac{ |E_i|^2}{ \beta |\partial E_i |}}
= \beta \frac{|\partial E _i| }{|E_i|^2} \Big ( 1 - \frac{ \Sigma _\infty ( E _i )}{\Sigma _\infty (E _i ) +  \frac{ |E_i|^2}{ \beta |\partial E_i |} } \Big ) \,.
$$ 
Therefore, the required estimate \eqref{t12} is satisfied provided 
$$\frac{ \Sigma _\infty ( E _i )}{\Sigma _\infty (E _i ) +  \frac{ |E_i|^2}{ \beta |\partial E_i |} }  \leq \e\, , $$
which holds since cells of Class 1 satisfy  $\beta |\partial E _i| \leq \e K _ \e \leq \e (\Sigma _\infty( E_i)) ^ {-1} |E_i|^2$. 
\medskip
\item[Class 2]: cells with $\beta |\partial E _i|  > \e K _\e$. For such cells, it holds  
\begin{equation}\label{t22}
\tau ^ {-1} ( E _i, \beta) \leq \beta \frac{|\partial E _i|}{|E_i| ^2} \leq \frac{4 ^5 \, \beta ^ 4 }{ \e ^ 5  K _ \e ^3 }\,. 
\end{equation}
Namely, using the upper bound \eqref{upperbo2} given by Proposition \ref{p:bounds2}, 
 the elementary estimates $|\partial E_i| \leq 4 d _i$, $|E _i| \geq \frac{1}{2} d _i w _i $,  and the fact that we are dealing with cells of Class 2 in $\conv (\e)$, 
 which satisfy in particular $4 d _i \beta \geq \e K _ \e$, 
 we get
$$\tau ^ {-1}  ( E _i,\beta)  \leq \beta \frac{|\partial E _i|}{|E_i| ^2} \leq \beta \frac{4 d _i}{\frac{1}{4} d _i ^2 w _i^2}   
= \frac{16 \beta  }{d_i w_i ^2}\leq \frac{16 \beta  }{ \e ^ 2 d_i ^3}  \leq \frac{16 \beta }{ \e ^ 2 \Big (  \frac{\e K _\e} {4 \beta} \Big ) ^3}  =  \frac{4 ^5 \, \beta ^ 4 }{ \e ^ 5  K _ \e ^3 }\,.
$$  
\end{itemize} 

Now, having at our disposal  the bounds \eqref{t12} and \eqref{t22} for cells of Class 1 and Class 2 respectively, the estimate in 
\eqref{claim2a2}
can be proved in the analogous way  as in Theorem \ref{t:honeycomb}. The idea is that 
the infimum which defines $m _{k, \e} (\Omega)$ has the same asymptotic behaviour if we restrict the  the sum  to indices $i \in \theta  _{k, \e}(\{E_i\})$   such that $E_i$ is of Class 1. 
It is enough to follow the proof of \eqref{claim2a},
with the obvious modifications in the scaling factors, and exploiting 
 \eqref{t12}-\eqref{t22} in place of \eqref{t1}-\eqref{t2}.

\medskip

{\it Step III: The lower bound \eqref{tesi12}
holds true}. 
On can repeat the same arguments used for Step III in the proof of Theorem \ref{t:honeycomb}, with the obvious modifications (in particular,  we exploit \eqref{aux12}, \eqref{claim12}, \eqref{claim2a2}, and  \eqref{behave2}).  

\bigskip
$\bullet$ {\bf Case $\beta <0$.} 
We address the reader to the proof of Theorem \ref{t:honeycomb} in case $\beta <0$, which runs exactly in the same way after suitably modifying the scaling factors. 
\qed
\section  {Proof of Corollary \ref{c:cormax}}\label{sec:proof3}
We give the proof only in case (i), as case (ii) is completely analogous. 

$\bullet$ {\bf Case $\beta >0$.} 
From the definition of $R_k (\Omega, \beta)$, we have $k R _k (\Omega, \beta) \geq r _k (\Omega, \beta)$, so that
$$
\limsup _{k\to + \infty} \frac{|\Omega|^ {1/2}}{ k ^ {1/2}} R _k (\Omega, \beta) 
\geq 
\limsup _{k\to + \infty} \frac{|\Omega|^ {1/2}}{ k ^ {3/2}} r _k (\Omega, \beta) =
\beta h ( H) 
\,,
$$
where in the last equality we have applied Theorem \ref{t:honeycomb}.

To prove the converse inequality, assume $|\Omega| = 1$, and  observe  that for every $\vps >0$ there exists $k_\vps$ such that for every $k\ge k_\vps$, the set $\Omega$ contains the convex $k$-cluster 
\begin{equation}\label{optimalcluster}
\big \{ \Big (\frac{1-\vps}{k}\Big )^\frac 12 {C(H)}  \, , \dots,  \Big (\frac{1-\vps}{k}\Big )^\frac 12 {C(H)} \big \}
\end{equation}
given by $k$ copies of the 
Cheeger set  $C (H)$ of the unit area regular hexagon, each one scaled so to have area $\frac{1-\vps}{k}$.  
We get
$$\begin{array} {ll} 
 \displaystyle \liminf _{k\to + \infty} \frac{R _k (\Omega, \beta) }{ k ^ {1/2} }  
& \displaystyle \leq
\liminf _{k\to + \infty}\frac{1 }{ k ^ {1/2} }   
\, \lambda _ 1 \Big ( \Big (\frac{1-\vps}{k}\Big )^\frac 12 C(H), \beta \Big )  \\ \noalign{\bigskip}
& \displaystyle
= \liminf _{k\to + \infty}\frac{1 }{ k ^ {1/2} }   \frac{k   }{1-\vps}
\, \lambda _ 1 \Big ( C(H) , \frac{\beta(1-\vps)^{1/2}}{k ^ {1/2}} \Big )  \\ \noalign{\bigskip}
& \displaystyle
= \liminf _{k\to + \infty}\frac{1 }{ k ^ {1/2} }   k   
 \frac{\beta}{ k ^ {1/2}}  \frac{1}{(1-\vps)^{1/2} }\frac {|\partial C (H)|} {|C (H)|} = \frac{1}{(1-\vps)^{1/2}} \beta h (H)  \,. \\ \noalign{\bigskip}
\end{array}
 $$
The parameter $\vps >0$ being arbitrary, we conclude the proof.

$\bullet$ {\bf Case $\beta <0$.} 
From the definition of $R_k (\Omega, \beta)$, we have in this case $k R _k (\Omega, \beta) \leq r _k (\Omega, \beta)$, so that
$$
\liminf _{k\to + \infty} \frac{|\Omega|^ {1/2}}{ k ^ {1/2}} R _k (\Omega, \beta) 
\leq 
\liminf _{k\to + \infty} \frac{|\Omega|^ {1/2}}{ k ^ {3/2}} r _k (\Omega, \beta) =
\beta h ( H) 
\,,
$$
where in the last equality we have applied Theorem \ref{t:honeycomb}.

To prove the converse inequality, we proceed as above. We assume without loss of generality that $|\Omega| = 1$, and we exploit the fact that 
asymptotically $\Omega$  contains   the convex $k$-cluster 
$\big \{ \Big (\frac{1-\vps}{k}\Big )^\frac 12 {C(H)}  \, , \dots,  \Big (\frac{1-\vps}{k}\Big )^\frac 12 {C(H)}  \big \}
$ to obtain
$$\liminf _{k\to + \infty} \frac{R _k (\Omega, \beta) }{ k ^ {1/2} }  
\geq
\liminf _{k\to + \infty}\frac{1 }{ k ^ {1/2} }   
k\, \lambda _ 1 \Big ( \frac{C(H)}{k^ {1/2}}, \beta \Big )  = \beta h (H)  \,.
 $$

\qed

\bigskip

\begin{remark}\label{finrem}
For $\beta >0$, a more direct proof of Theorem \ref{t:honeycomb} (and consequently of Corollary \ref{c:cormax} (i)) can be performed by using the following lower bound in place of \eqref{lowerbo}:
\begin{equation}\label{newlowerbo}
\lambda _ 1 (\Omega, \beta) \geq \beta h (\Omega ) - \beta ^ 2 \,.
\end{equation}
 We point out that this inequality, which must be attributed to Bossel \cite{Bossel}, holds true without any assumption on $\Omega$. It is in general quite rough (for instance, the right hand side can have negative sign  for some $\Omega$ and $\beta$), but it becomes useful as soon as the Cheeger constant $h(\Omega)$ becomes large, which is typically the case in a partition with a large number of cells. For the sake of completeness, we enclose a short independent proof. 
Recall that the functional formulation of the Cheeger constant  reads
\begin{equation}\label{minCheg}
h (\Omega) = \inf _{v \in BV (\R ^ 2 ) \setminus \{ 0 \}, v=0 \mbox{ on }  \R ^ 2 \setminus \Omega} \frac{ |Dv| (\R ^ 2)} {\int _\Omega |v| } \,.
\end{equation}
Then, letting $u$ be a first Robin eigenfunction, extended to $0$ out of $\Omega$, and taking $v = u ^ 2$ in the  minimization problem \eqref{minCheg}, gives
$$\begin{array}{ll}
\beta h (\Omega) & \displaystyle \leq  \frac{ \beta \int_\Omega |\nabla ( u ^ 2) | + \beta \int _{\partial \Omega} u ^ 2} {\int _\Omega u ^ 2 }  =  \frac{\beta \int_\Omega 2 | u \nabla  u  | +\beta  \int _{\partial \Omega} u ^ 2} {\int _\Omega u ^ 2 }  
\\
\noalign{\medskip}
& \displaystyle \leq  \frac{ \int_\Omega  |  \nabla  u  |^2  + \beta ^ 2 \int _\Omega u ^ 2 +\beta  \int _{\partial \Omega} u ^ 2} {\int _\Omega u ^ 2 }   =  
\lambda _ 1 (\Omega, \beta) + \beta ^ 2 \,, 
 \end{array}
$$
which is exactly \eqref{newlowerbo}. Now, 
by exploiting such estimate, it is immediate to obtain the required inequality \eqref{tesi1} in the proof of Theorem \ref{t:honeycomb}. 

We emphasize that this more direct approach does not work for the Robin torsional ridigity,  which still requires the finer inequality on convex sets. Indeed,  Bossel's approach seem to be unadaptable to the case of torsional rigidity; as well, 
by arguing as above one arrives at the lower bound 
\begin{equation}\label{newlowerbo2}
\tau ^ {-1} (\Omega, \beta) \geq \frac {\beta}{|\Omega|} h (\Omega ) - \frac{\beta ^ 2}{|\Omega|} \,, 
\end{equation}
which does not allow obtain the required inequality \eqref{tesi12} in the proof of Theorem \ref{t:honeycomb2}. 
Instead, \eqref{newlowerbo2} may be used to prove the statements analougue to
Theorem \ref{t:honeycomb2} and Corollary \ref{c:cormax} (ii) for the more unusual functional $\lambda (\Omega, \beta) : = |\Omega| \tau ^ {-1} (\Omega, \beta)$.

\end{remark}

\bigskip
\bibliographystyle{mybst}
\bibliography{References}

\end{document}